\documentclass[12pt]{article}
\usepackage{latexsym}
\usepackage{amsfonts}
\usepackage{amsmath}
\usepackage{amssymb}

\makeatletter

\ifx\reset@font\undefined
             \let\reset@font=\relax
             \fi
\def\section{\@startsection {section}{1}{\z@}{-3.5ex plus-1ex minus
               -.2ex}{2.3ex plus.2ex}{\reset@font\large\bf}}
\def\subsection{\@startsection{subsection}{2}{\z@}{-3.25ex plus-1ex
               minus-.2ex}{-.1em}{\reset@font\large\bf}}
\def\subsubsection{\@startsection{subsubsection}{3}{\z@}{-3.25ex plus
            -1ex minus-.2ex}{-.1em}{\reset@font\normalsize\bf}}
\makeatother

\title{On the nontrivial projection problem }
 \author{Stanislaw J. Szarek
 \thanks{Supported in part  by grants from the
 National Science Foundation (U.S.A.).}
                    \and  Nicole Tomczak-Jaegermann
 \thanks{This author holds the Canada
       Research Chair in Geometric  Analysis.}
 }

 \newcommand\address{\noindent\leavevmode%
 Department of Mathematics\\
 Case Western Reserve University\\
 Cleveland, Ohio 44106-7058, U.S.A.\\
 {\it and}\\
 Equipe d'Analyse Fonctionnelle, BP 186\\
 Universit\'e Pierre et Marie Curie\\
 75252 Paris, France\\
 {\small\tt%
 szarek@cwru.edu}\\[.5cm]
 \noindent
 Department of Mathematical and Statistical Sciences,\\
 University of Alberta,\\
 Edmonton, Alberta, Canada T6G 2G1,\\
 {\small\tt%
              nicole@ellpspace.math.ualberta.ca} }

\date{}

\newtheorem{fact}{Fact}   
\newtheorem{thm}[fact]{Theorem}
\newtheorem{prop}[fact]{Proposition}

\newbox\nrmbox
\setbox\nrmbox=\hbox{$\Vert$}
\def\nrmrule{\vrule height\ht\nrmbox depth1.2\dp\nrmbox}
\setbox\nrmbox=%
             \hbox{\kern0.15em\nrmrule\kern0.15em\nrmrule\kern0.15em%
                       \nrmrule\kern0.15em}
\newcommand{\Snorm}[1]%
             {\copy\nrmbox#1\copy\nrmbox\kern-0.03em%
                   \lower.4ex\hbox{}}

\newcommand{\R}{\mathbb{R}}
\newcommand{\Rn}{\R^n}
\newcommand{\N}{\mathbb{N}}
\newcommand{\E}{\mathbb{E}}

\newcommand{\ep}{\varepsilon}

\newcommand{\rank}{\mathop{\rm rank\,}}
\newcommand{\codim}{\mathop{\rm codim\,}}

\begin{document}

\maketitle

{\abstract{\small
The {\em Nontrivial Projection Problem} asks whether 
every finite-dimensional normed space admits a well-bounded
projection of non-trivial rank and corank or, equivalently, whether 
every centrally symmetric convex body (of arbitrary dimension) is approximately 
affinely equivalent to a direct product of two bodies of non-trivial dimension. 
We show that this is true ``up to a logarithmic factor."}}

\section{Introduction and the main results }
\label{intro}

A series of well-known open problems in asymptotic theory of
normed spaces is concerned with the existence, in any
finite-dimensional normed space (of dimension greater than one), of well-bounded
projections of non-trivial rank and corank. One possible 
formulation is as follows

\medskip
\noindent {{\bf The Nontrivial Projection Problem}} \ \ {\sl
  Do there exist $C\geq 1$ and a sequence $k_n \to \infty$ such
  that for every $n$-dimensional normed space $X$ there is
  a projection $P$ on $X$ with \\
 {\rm (i)}  \,$\|P\|\leq C$  \\
 {\rm  (ii)} $\min\{{\rm rank} \,P, {\rm rank} (I-P)\} \geq k_n$ ?}

\medskip \noindent 
Versions of this question were explicitly posed in
ICM talks by Pisier (\cite{PiICM}, 1983) and Milman
(\cite{MiICM86}, 1986).  In geometric terms,
the problem asks whether, for $n\geq 2$,
an arbitrary $n$-dimensional centrally symmetric convex body 
 is approximately (``up to a constant $C$")
affinely equivalent to a direct product of two bodies 
 whose dimensions are at least $k_n$.

\medskip
To put the problem in a perspective, for a subspace $E $ of a
Banach space $ X$ denote
$$
\lambda(E, X) := \inf \{ \|P\| : P \mbox{ is a projection from
} X \mbox{ onto } E \} .
$$
 We then have  (Kadets and Snobar \cite{KS}, 1971)
$$
\dim E = k  \Rightarrow  \lambda(E, X) \leq \sqrt{k}
$$
or, more precisely (K\"onig and Tomczak-Jaegermann
\cite{KT}, 1990),
$$ 
\lambda(E, X) \leq \sqrt{k} - c/\sqrt{k}
$$
for all $k>1$ and some universal (and explicit) numerical constant $c>0$.

\medskip The estimates above hold for all subspaces, and
sometimes can not be substantially improved.  First, there
is the remarkable {\em infinite dimensional} example of
Pisier (\cite{PiActa}, 1983):

\medskip \noindent {\bf Pisier's space} {\sl There exists a
  Banach space $X$ and $c>0$ such that for any finite rank
  projection $P$ on $X$ one has $\|P\|\geq c \,\sqrt{{\rm
      rank} P}$. }

\bigskip \noindent 
Next, it follows from the work on the finite
dimensional basis problem (Gluskin \cite{Glu2}, 1981; Szarek
\cite{Sz83}, 1983) that, in general, we may not be able to
find projections on $X$ whose rank and corank are of the
same order as $\dim X$ and whose norm is $o(\sqrt{\dim X})$,
and that the statement from the problem can not hold with
$k_n$ substantially larger than $\sqrt{n}$ (more precisely,
with $k_n \gg \sqrt{n \log n}$).

\medskip However, all these results do not exclude a positive
answer to the following (sample) question.

\medskip \noindent {\bf A generalized Auerbach system } {\sl
  Does there exist $C\geq 1$ such that for any $n \in \N$, for
  any $n$-dimensional normed space $X$, and for any integer
  $m$ with $\sqrt{n} < m \leq n$, the space $X$ can be split
  into a direct sum of $m$ subspaces $E_1, \ldots , E_m$ of
  approximately equal dimensions, and such that if $P_j$ is
  the projection onto $E_j$ that annihilates all $E_i$'s with
  $i \neq j$, then $\max_{1\leq j \leq m} \|P_j\|~\leq~C$?}

\medskip \noindent A positive answer would of course imply a
positive solution to the nontrivial projection problem.  The
{\em classical } Auerbach lemma asserts that the answer is
``yes, with $C=1$," if $m=n$.

\medskip In the positive direction, it has been known for
quite a while that in some cases bounds on the norm sharper
than $\sqrt{\min\{{\rm rank} \, P, {\rm rank} (I-P)\}} $ 
can be obtained, primarily via
arguments based on $K$-convexity (Figiel and
Tomczak-Jaegermann \cite{FT}, 1979; Pisier \cite{PiAENS},
1980).  Based on that point of view and on the arguments and
results from \cite{PiActa}, Pisier posed (\cite{PiCle,
  PiActa}) modified variants of the nontrivial projection
problem.  One possible formulation is the following version
of the uniformly complemented $\ell_p^n$ conjecture of
Lindenstrauss \cite{Li}.

\bigskip\noindent {{\bf The modified problem}} {\sl Given a
  sequence $(X_n)$ of finite dimensional normed spaces with
  $\dim X_n \to \infty$, does there exist $p \in \{1, 2,
  \infty\}$, a constant $C\geq 1$ and sequences $m_k\to \infty$
  and $n_k\to \infty$ such that $X_{n_k}$ contains a
  subspace which is $C$-complemented and $C$-isomorphic to
  $\ell_p^{m_k}$? }
    
\medskip \noindent 
It is worthwhile to note that, up to the precise
value of the constant, the conditions on the subspace can be
conveniently rephrased as ``{\em the identity on
  $\ell_p^{k}$ $C$-factors through $X_{n_k}$.}" 
  (We refer to the next section for definitions of concepts that
  may be unfamiliar to a non-specialist reader.)

\medskip An affirmative answer to the modified problem would
follow from an affirmative answer to the following (see,
e.g., \cite{PiKconvex, PiCle})

\medskip \noindent {{\bf The cotype-cotype conjecture}} {\sl
  If a Banach space $X$ has (an appropriate) approximation
  property and if both $X$ and its dual $X^*$ have
  nontrivial cotype, then $X$ is $K$-convex.}

\medskip \noindent 
Pisier's example mentioned earlier shows that {\em some} approximation 
hypothesis is necessary.
The setting that is of interest to us is finite dimensional, with dimension-free
estimates on the parameters involved, and so the issues related to 
approximation properties will not enter the discussion.

\medskip  In the present paper we shall prove the following 
result in the  direction of the nontrivial projection problem.

\begin{thm} \label{main} There exist $C >0$ and a sequence
  $k_n \to \infty$ such that, for every $n\geq 2$ and for every $n$-dimensional
  normed space $X$, there is
  a projection $P$ on $X$ with \\
  {\rm (i) }  $\|P\|\leq C (1+\log k_n)^2$  \\
  {\rm (ii)} $\min\{{\rm rank} \,P, {\rm rank} (I-P)\} \geq
  k_n$.\\ 
  Moreover, the range of the projection $P$ is
  $C$-isomorphic to an $\ell_p$-space for some $p \in \{1,
  2, \infty\}$.
\end{thm} 

\noindent
{\bf Remarks} (a) The argument shows that one can 
choose $(k_n)$ to grow  as (roughly) $\exp(\sqrt{\log n})$. \\ 
(b) A slightly weaker but more compact statement than the
assertion of the Theorem is
``{\em the identity on  $\ell_p^{k_n}$ can be $C (1+\log
  k_n)^2$-factored through $X$}."\\ 
(c) Already in this last form, the assertion is nearly
optimal (even  for the class of $\ell_q^n$ spaces), except
for the exact values of the powers of $\log k_n$ in
(i). Similarly, $k_n$ can not be substantially larger than
the quantity given in Remark (a) above.  This is explained in
section \ref{opt}; see also the remark following the proof
of the Theorem.

\medskip The proof of the Theorem is based on a dichotomy
which yields either

\smallskip \noindent
(1) a reasonably complemented copy of $\ell_2^{k_n}$ via an
argument based on $K$-convexity and the $\ell$-ellipsoid
(essentially as in \cite{FT}) {\em or}\\
(2) a good copy of $\ell_\infty^{k_n}$ in $X$ or in $X^*$,
necessarily well-complemented (the latter implies existence
of a well-complemented copy of $\ell_1^{k_n}$ in $X$).  This
part is based on a result of Alon-Milman with a refinement
due to Talagrand, on restricted invertibility results in the
spirit of Bourgain-Tzafriri, on a blocking argument due to
James, and on various tricks of the trade developed over the
last 25 years.

\section{Notation and preliminaries} \label{prelim}

We use the standard notation from the Banach space theory. 
In particular, we denote Banach (or normed) spaces by $X$, $Y$ etc.,
 and by $B_X$, $B_Y$ \ldots their (closed) unit balls.  
 An operator means a bounded linear operator. 
 For an operator $T: X \to Y$, its
operator norm is denoted by $\|T:X \to Y\|$ or just by $\|T\|$.
For isomorphic Banach spaces $X$ and $Y$, their Banach-Mazur
distance is defined by $d(X, Y)= \inf \|T\|\, \|T^{-1}\|$,
where the infimum is taken over all isomorphisms $T$ from
$X$ onto $Y$; we say that $X$ is $\lambda$-isomorphic to $Y$ if
$d(X, Y)\le \lambda$. A subspace $F$ of $X$ is
$\lambda$-complemented if there exists a projection from $X$ onto
$F$ of norm less than or equal to $\lambda$.

For finite-dimensional normed spaces, the essentially
equivalent language of symmetric convex bodies is natural and
often very useful. (By a symmetric convex body $K
\subset \R^n$ we will mean a convex compact set
with non-empty interior which is centrally
symmetric with respect to the origin.) 
By $\|~\cdot~\|_K$ we denote the gauge of $K$; 
then $X = (\R^n, \|\cdot\|_K)$ is an
$n$-dimensional normed space such that $K = B_X$.  Any $n$-dimensional
normed space can be represented in such a form in many
different (although isometric) ways.  
If $K_1 \subset \R^{n_1}, K_2  \subset \R^{n_2}$ are
symmetric convex bodies and $X_1, X_2$ are the
corresponding normed spaces, for an operator $T:\R^{n_1} \to
\R^{n_2}$ the operator norm $\|T: X_1\to X_2\|$ will be also denoted
by $\|T: K_1\to K_2\|$ or (for example) by $\|T: K_1\to X_2\|$.

By $|\cdot|$ we denote the Euclidean norm on $\R^n$ and we
use the representation $\ell_2^n = (\R^n, |\cdot|)$. The
Euclidean ball in $\R^n$ and the inner product are denoted
by $B_2^n$ and $\langle \cdot, \cdot \rangle$.  For a
subspace $E \subset \R^n$,  we denote by $P_E$ the orthogonal
projection on $E$.  The polar body $K^\circ$ is defined by
$K^\circ := \{ x \in \R^n \mid |\langle x, y\rangle|\le 1 \
\mbox{\rm for all}\ y\in K\}$.  As is well known, the normed space
$(\R^n, \|\cdot\|_{K^\circ}) $ can be canonically identified with the dual
space $(\R^n, \|\cdot\|_{K})^*$.

We now recall the following less standard concept which will be
useful further on.  Given a normed space $Y$ and a linear
operator $S : \ell_2^n \to Y$, the $\ell$-norm of $S$ is
defined via 
$ \ell(S) := \left(\int_{\R^n} \|Sx\|^2 \,
  d\mu_n\right)^{1/2}$, where $\mu_n$ is the standard
Gaussian measure on $\R^n$. In other words,
$$
\ell(S) = \left(\E  \|\sum_{i=1}^n g_i Sv_i\|^2\right)^{1/2},
$$
where $(v_i)$ is an arbitrary orthonormal basis of
$\ell_2^n$ and $(g_i)$ -- an i.i.d. sequence of $N(0,1)$
Gaussian random variables
($\E $ stands for the expected value). It is well-known and easy to verify that 
the $\ell$-norm satisfies $\|S\|\le \ell(S)$ and it 
has the ideal property  $\ell(SA) \leq \ell(S)\|A:
\ell_2^n \to \ell_2^n \|$. We refer the
reader to \cite{TJ} or \cite{PiBook} for more details.

For a symmetric convex body  $K\subset \R^n$ we set 
$$
\ell(K) = (\E \|g\|_K^2)^{1/2},
$$
where $g = (g_1, \ldots, g_n) \in \R^n$ is the standard
Gaussian vector. (In other words, $\ell(K) = \ell(J)$ where
$J:\ell_2^n \to K$ is the formal identity operator.)  It is
known
 that
one may find a linear image $\tilde K = u K$ (with $u:\R^n
\to \R^n$ one-to-one and onto), called by some authors the
$\ell$-position of $K$, which in particular satisfies
\begin{equation}
  \label{ell-ellips}
\ell(\tilde K )= \ell((\tilde K)^\circ)
 \leq C\sqrt{n (1+\log{n})}\, ,
\end{equation}
where $C$ is a universal constant.  Clearly, 
the normed space induced by $\tilde K$ 
is isometric to the one associated with $K$.

The inequality (\ref{ell-ellips}) lies at the core of our arguments. 
Is is obtained by combining results of \cite{FT} and \cite{PiKconvex}, in
particular by using deep connections with $K$-convexity; 
this is where the $\ell$-position/$\ell$-ellipsoid come in.

As in inequality (\ref{ell-ellips}) and earlier in the introduction, the symbols 
$c,C,c', C_1$ etc. will stand in what follows for {\em universal} positive constants, 
independent of the particular instance of the problem  that is being considered
(most notably independent of the dimension).  However, the same symbol may  
represent different numerical values in different parts of the paper.

\section{Proof of the Theorem} The argument will be based on
two Propositions corresponding to the two alternatives of
the dichotomy mentioned in the Introduction.

\begin{prop} \label{dvor}
Let $K_1, K_2 \subset \Rn$  be symmetric convex bodies such that  $K_1 \subset \alpha B_2^n$ 
and $K_2 \supset \beta^{-1} B_2^n$ and let $\sqrt{m} \leq c \,\min\{\ell(K_1^\circ))/\alpha, \ell(K_2)/\beta\}$.
Then, for most of  subspaces $F$ of $\Rn$ of dimension $m$ (in the sense of the Haar measure on the corresponding Grassmannian), \\
{\rm (i) } $\|P_F : K_1 \to K_2\| \leq C\,\ell(K_1^\circ))\ell(K_2)/n$\\
{\rm (ii)} $\exists \; r>0$ such that $r(B_2^n \cap F) \subset K_2 \cap F \subset Cr(B_2^n \cap F)$.
\end{prop}

Here is a sketch of the proof based on Milman's version of
Dvoretzky theorem (\cite{M1} or \cite{MS}, Chapter
4). First, if $m \leq (c \,\ell(K_2)/\beta)^2$, then, for
most of subspaces $F$ of dimension $m$, the section $K_2
\cap F$ is approximately a Euclidean ball of radius $r =
\sqrt{n}/\ell(K_2)$, which yields (ii).  Dually, if $m \leq
(c\,\ell(K_1^\circ)/\alpha)^2$, then $P_F K_1$ is -- again,
for most $F$'s -- approximately a Euclidean ball of radius
$R = \ell(K_1^\circ)/\sqrt{n}$.  If $F$ is such that both of
the above hold, then $\|P_F : K_1 \to K_2\|$ is
approximately $R/r$, whence (i) follows.

We point out that most authors use in similar arguments {\em
  spherical} rather that {\em Gaussian} averages; this 
is why our formulae involve $\sqrt{n}$ factors that are
absent, e.g., in \cite{MS}.

\medskip 
The second technical result that we need is the following.
\begin{prop} \label{ell1}
Let $K_0 = B_Y \subset \Rn$  be such that 
for any subspace $E \subset \Rn$ with 
$\codim E < k$ we have $\|P_E : Y \to \ell_2^n\| \geq a$. Then there exists a subspace 
$Z$ of $Y$ such that $d(Z,\ell_1^m) \leq C$ with 
$m:= \dim Z \geq  c \, k^{1/\gamma}$, 
where $\gamma = 2\log_2({32\,\ell(K_0^\circ)}/{a })$ 
and a projection $Q : Y\to Z$ with $\|Q\| \leq C$.
\end{prop}

We postpone the proof of Proposition \ref{ell1} until the next section and direct our
attention to Theorem \ref{main}.  The argument will split naturally into three parts 
corresponding to different choices of  $p \in \{1,2,\infty\}$, which will in turn depend
 on the values of 
certain parameters related to the geometry of $X$.  

Set $K=B_X$ and $k = \lceil n/4 \rceil$. 
Let $a  \in [1, \sqrt{n}]$ ($a$ will be later specified to be roughly $\sqrt{n}/ \exp{\sqrt{\log n}}$ ).

Assume now that for every subspace $E \subset \Rn$ with 
$\codim E < k$ we have $\|P_E : X \to \ell_2^n\| \geq a$.  Accordingly, Proposition 
\ref{ell1} applies for $Y=X$, yielding a well-complemented 
$m$-dimensional subspace of $X$, well-isomorphic to $\ell_1^m$, with 
$m \geq c\,k^{1/\gamma} \geq c'n^{1/\gamma}$, where 
$\gamma = 2\log_2({32\,\ell(K^\circ)}/{a })$. 

\smallskip
Similarly, if,  for every subspace $E \subset \Rn$ with 
$\codim E < k$ the estimate $\|P_E : X^* \to \ell_2^n\| \geq a$ holds, then the same argument 
produces a well-complemented subspace of $X^*$ well-isomorphic to $\ell_1^m$.
By duality, this yields  a well-complemented subspace of $X$ well-isomorphic to 
$\ell_\infty^m$, with the bound for $m$ involving now 
$\gamma = 2\log_2({32\,\ell(K)}/{a })$. 

\medskip
If neither of these conditions is satisfied, then there exist subspaces $E_1, E_2 
\subset \Rn$ of codimension $< n/4$ such that the appropriate norms of projections 
$P_{E_1}, P_{E_2}$ do not exceed  $a$, and so also $\|P_H : X \to \ell_2^n\| \leq a$ and 
$\|P_H : X^* \to \ell_2^n\| \leq a$,   where $H=E_1\cap E_2$.  In geometric terms, 
this is equivalent to the inclusions 
$$
P_H K \subset a B_2^n, \ \ \ K \cap H  \supset a^{-1}  (B_2^n \cap H),
$$
the latter of which is the dual reformulation of $P_H K^\circ \subset a B_2^n$ .
We are thus in a position to apply Proposition  \ref{dvor} with $K_1 = P_H K$, 
$K_2 = K \cap H$, $\alpha = \beta = a$ and $H$ playing the role of $\Rn$. 
(Note that $\dim H > n/2$.)  This yields existence of a $C$-Euclidean section 
$K_2 \cap F = K \cap F$, whose dimension $m$ is of order 
$\left(\min\{\ell((P_H K)^\circ), \ell(K \cap H)\}\right)^2/a^2$. 
Moreover, $P_F K = P_F(P_H K) \subset \lambda (K \cap H) \subset \lambda K$,
where $\lambda \leq C \ell((P_H K)^\circ)\ell(K \cap H)/n$. In other 
words, $F$ is a $\lambda$-complemented $C$-Euclidean subspace of $X$.

\smallskip It remains to collect estimates on ranks and
norms of the projections and choose an optimal value for
$a$. This will also require choosing an appropriate
representation of $X$ on $\Rn$, namely the $\ell$-position,
so that condition  (\ref{ell-ellips}) of Section \ref{prelim} is satisfied. In
particular, we will have
 $$
 \ell(K\cap H) \leq  \ell(K) = \sqrt{n \kappa} , 
 $$
 $$
 \ell((P_HK)^\circ)=   \ell(K^\circ \cap H) \leq  \ell(K^\circ) = \sqrt{n \kappa} ,
 $$
 where $\kappa \leq C (1+\log n)$. 
On the other hand, 
$$
\frac n 2 < \dim H \leq  \ell((P_HK)^\circ)  \ell(P_HK) \leq  \ell(K^\circ \cap H) \ell(K \cap H),
$$
and so we also have lower estimates 
$$
 \ell(K\cap H) \geq  \frac 12   \sqrt{\frac {n}{\kappa}}, \ \ 
  \ell((P_HK)^\circ)=   \ell(K^\circ \cap H) \geq \frac 12   \sqrt{\frac {n}{\kappa}} .
$$

The lower bounds for dimensions of subspaces become 
$$
\frac {c' n}{\kappa a^2}, \ \ c'n^{1/{2\log_2({32\sqrt{n\kappa}/a})}} ,
$$
for $p=2$ and $p=1$ or $ \infty$ respectively.  Choosing $a = \sqrt{n}/ \exp{\sqrt{\log n}}$
and remembering the upper bound on $\kappa$ we easily check that both of these 
quantities are $\geq k_n:=c \exp{(\frac 12 \sqrt{\log n}}) $. On the other hand, 
the upper bound on the norm of projection in the case $p=2$  is clearly 
$C_1\kappa \leq C_2 (1+ \log n) \leq C_3 (1+\log k_n)^2$, which concludes the proof of the 
Theorem together with the bound on $k_n$ given in Remark (a).

\medskip
Our final comment concerns optimality of the estimate for the norms of the projections 
in terms of their rank. By choosing differently the threshold value $a$, we can 
increase the dimension of  the $C$-Euclidean subspace $F$, while keeping the 
norm of $P_F$ bounded by $C\log n$.    This way we can assure 
that, in all cases, the norm of the projection $P$ is $\leq C' \log (\rank P)  \log \log (\rank P) $.
The price we pay is a decrease in the 
dimensions of the $\ell_1^m$ or $\ell_\infty^m$ subspaces, and the common lower bound 
for ranks of projections is only a power of $\log n$ instead of  $\exp{(c\sqrt{\log n})}$.
(The power of $\log n$ can be chosen arbitrarily, at the cost of increasing the constant $C'$.)

\section{Proof of Proposition \ref{ell1}}
%
We start by defining (by induction) two sequences $x_1, x_2, \ldots, x_k$ and $y_1,
y_2, \ldots, y_k$ with certain extremal properties.  First, let $x_1
=y_1 \in K_0$ be such that $|x_1|=a_1 := \max _{x \in K_0}|x|$.
For consistence with future notation set $F_1 = \R^n$.
  Next, suppose that $1<j\leq k$
and that $x_i, y_i$ for $i<j$ have already been defined.  Set  $F_j
:= [x_1, x_2, \ldots, x_{j-1}]^\perp$ and choose $y_{j} \in K_0$ so that
$|P_{F_j} y_{j} | = \|P_{F_j} : K_0 \to \ell_2^n\| =: a_{j}$. 
Set $x_{j} = P_{F_j} y_{j} $; then the sequence $(x_j)$ is orthogonal
with $|x_j|=a_j$.  Finally, define an
orthonormal sequence $(u_j)$ by  $u_j:=x_j/a_j$, $j=1, 2, \ldots, k$. 
Note that, by hypothesis and construction, $a_1 \geq a_2\geq \ldots a_k\geq a$.

Pick an interval $I \subset \{1, \ldots, k\}$ with $|I| \ge k/(1+ \log_2
\frac{a_1}{a_k} )$ such that $a_i \le 2 a_{i'}$ for all $i, i' \in I$.
Set $F := [x_i]_{i \in I} $, and let $a' = a_{\min I}$.  In the sequel
we will analyze the convex set $\tilde{K_0} :=P_F
K_0 $-- viewed as a convex body in $F$ -- and sequences $\tilde{y}_j
:= P_F y_j$. 
By construction, all
$\tilde{y}_j$'s are elements of $\tilde{K_0}$.  Moreover, since $\|P_F: K_0
\to B_2^n\| \le \|P_{F_{\min I}}: K_0 \to B_2^n\| = a'$, it follows that for all
$w \in F$,
\begin{equation} \label{ell2vsK}
|w| \leq  a' \; \|w\|_{\tilde{K_0}}
\end{equation}
and, in particular,  $|\tilde{y}_j| \leq a'$ for $j \in I$.
On the other hand, since 
$P_{F_j} u_j = P_{F} {u_j} =  u_j $ 
for $j \in I$, it follows that for such $j$
$$\langle \tilde{y}_j, u_j\rangle =
\langle P_{F} {y}_j, u_j\rangle = \langle{y}_j, u_j\rangle = 
\langle P_{F_j} {y}_j, u_j\rangle = 
\langle {x}_j, u_j\rangle = 
a_j\ge a'/2 .$$
Accordingly, we are in a position to apply Bourgain-Tzafriri
restricted invertibility principle (\cite{BT}) in the form presented
in Lemma B in \cite{BS} to conclude that there exists a set $\sigma
\subset I $ such that $s := |\sigma| \geq ck/(1 + \log_2
\frac{a_1}{a_k} )$ and verifying, for any sequence of scalars
$(t_j)_{j \in \sigma}$,
$$
\left|\sum_{j \in \sigma} t_j \tilde{y}_j \right| \geq 
\frac{a'}8 \left(\sum_{j \in \sigma}|t_j|^2\right)^{1/2}.  
$$

To reduce the clutter of subscripts, 
we will assume that $\sigma = \{1,2,\ldots, s\}$.
Let $(z_j)_{j=1}^s$ be the sequence in $[\tilde{y}_1, \tilde{y}_2, \ldots, \tilde{y}_s]$ that 
is biorthogonal to  $(\tilde{y}_j)_{j=1}^s$, then 
\begin{equation} \label{upperell2}
\left|\sum_{j=1}^s t_j z_j \right| \leq 
\frac8{a'} \left(\sum_{j=1}^s|t_j|^2\right)^{1/2}
\end{equation}
for any sequence of scalars $(t_j)$.   

Next, consider two polar bodies: $\tilde{K^\circ_0}$, the
polar of $\tilde K_0 $ {\em inside}  $F$, and ${K^\circ_0}$, the polar of $ K_0 $
(in $\R^n$).  Since $\tilde K_0$ is an orthogonal projection of $K_0$, i.e.,
$\tilde K_0 = P_F K_0$, it follows that $\tilde K_0^\circ$ is a section of $K_0^\circ$,
namely $\tilde K_0^\circ = K_0^\circ \cap F$.  Thus, given that $\|
\tilde{y}_j\|_{\tilde{K_0}} \leq 1$ and $z_j \in F$, it follows that
$\|z_j\|_{K_0^\circ} = \|z_j\|_{\tilde K_0^\circ} \geq 1$ for $1 \le j
\le s$.

Consider now the quantity $M:=\left(\E \|\sum_{j=1}^s
  g_iz_i\|_{K_0^\circ}^2\right)^{1/2}$
and define a linear map $T : \ell_2^n \to \ell_2^n$ by
$Te_j=z_j$ for $j=1,2,\ldots, s$ and $Te_j=0$ for $j>s$; it
then follows from (\ref{upperell2}) that $\|T\| \leq 8/a'$.
Accordingly, denoting by $J$ the identity map considered as an
operator from $\ell_2^n$ to $Y^* = (\Rn, \|\cdot
\|_{K_0^\circ})$,
and using the definition and properties of the $\ell$-norm discussed in
Section \ref{prelim}, we conclude that
\begin{eqnarray}
M&= &\left(\E  \|\sum_{j=1}^n g_i 
  Te_i\|_{{K^\circ_0}}^2\right)^{1/2} =\ell(J  T) \nonumber
\\  
&\leq &\ell(J) \| T : \ell_2^n \to \ell_2^n\| \leq  \ell(J)
\frac 
8{a'} =  \frac 8{a'}  \ell({{K^\circ_0}}) .  \label{em}
\end{eqnarray}
We now want to appeal to \cite{Tal} to extract from $(z_i)$
a subsequence resembling an $\ell_\infty$ basis. To this
end, we need to consider the modified average $M_1:= \E
\|\sum_{j=1}^s \ep_iz_i\|_{{K^\circ_0}}$, where $(\ep_i)$ is
an i.i.d. sequence of Bernoulli random variables. As is
well-known, $M_1 \leq \left(\E \| \sum_{j=1}^s \ep_i z_i
  \|_{{K^\circ_0}}^2 \right)^{1/2} \leq M$, which combined
with (\ref{em}) yields
\begin{equation} \label{em1}
M_1= \E  \|\sum_{j=1}^s \ep_iz_i\|_{{K}^\circ_0} \leq  \frac
8{a'}  \ell({{K^\circ_0}}) \leq \frac
8{a}  \ell({{K^\circ_0}}). 
\end{equation}
Another quantity that is needed to appeal to \cite{Tal} is $w := \max 
\|\sum_{j=1}^s \ep_iz_i\|_{{K^\circ_0}}$ (i.e., the maximum over all
choices of $\ep_i=\pm 1$). Dualizing estimate (\ref{ell2vsK}) and
using (\ref{upperell2}), we obtain, for all such $(\ep_i)$,
$$
\left\|\sum_{j=1}^s \ep_iz_i\right\|_{{K^\circ_0}}
= \left\|\sum_{j=1}^s \ep_iz_i\right\|_{\tilde{K^\circ_0}} \leq 
a' \left|\sum_{j=1}^s \ep_iz_i \right|
\leq 8 \;  s^{1/2}.
$$
We are now ready to use the following result from \cite{Tal}.

\begin{fact} \label{talafact}
Let $(z_i)_{i=1}^s$ be a sequence in a normed space.  Set
  $M_1 = \E \|\sum_{j=1}^s \ep_iz_i\|$ and $w= \max \|\sum_{j=1}^s
  \ep_iz_i\|$.  Then there exists a subset $\tau \subset \{1,2,\ldots
  , s\}$ with $|\tau| \geq s M_1/2w$ such that, for any scalars
  $(t_i)$,
$$
\|\sum_{i\in \tau} t_i z_i \| \leq 4 M_1 \max_{i\in \tau} |t_i | .
$$
\end{fact} 
Specified to our context, the Fact yields $\tau$ with $|\tau| 
\geq c'\left(\frac {k}  {1 + \log  \frac{a_1}{a_k}}\right)^{1/2} M_1$.

\smallskip
The next step is a well-known blocking argument due to R. C. James.

\begin{fact} \label{james} Let $v_1, v_2,\ldots,v_{m^2}$ be elements
  of a normed space $V$ 
  with $\|v_j\| \geq 1$ for all $j$ verifying, for some $\beta \geq
  1$, $\| \sum_{j=1}^{m^2} t_j v_j \| \leq \beta \max_j |t_j|$ for all
  sequences of scalars $(t_j)$. Then there exist $v_1',
  v_2',\ldots,v_{m}'  \in V$
   with $\|v_i'\| \geq 1$ for all $i$ such
  that $\| \sum_{i=1}^{m} t_i v_i' \| \leq \beta^{1/2} \max_i |t_i\|$
  for all sequences of scalars $(t_i)$.
\end{fact}
The  proof of Fact \ref{james} is based on the following dichotomy. If
there is a subset $\sigma \subset \{1,2,\ldots,n^2\}$ with $|\sigma|=m$,
for which $\| \sum_{i\in \sigma} \pm v_i \| \leq \beta^{1/2} $ for all choices of 
signs, then the collection $\{ v_i \, : \, {i\in \sigma} \}$ works.
If not, then for each such $\sigma$ there is $v_\sigma = \beta^{-1/2} \sum_{i\in \sigma} \pm v_i$
with $\|v_\sigma\| > 1$; partitioning
the set $\{1,2,\ldots,n^2\}$ into subsets 
$\sigma_1, \sigma_2,\ldots,\sigma_m$ with $|\sigma_j|=m$ for all $j$
we are led to a collection $\{ v_{\sigma_1}, v_{\sigma_2}, \ldots, v_{\sigma_n}\}$
which has the required property.

\medskip 
The procedure implicit in Fact \ref{james} can clearly be
iterated. Applying it $d=\lfloor \log_2 \log_2 (4M_1) \rfloor$ times to
our sequence $(z_i)_{i\in \tau}$ we are led to $z_1', z_2',\ldots ,
z_l'$ such that $\|z_i'\|_{K_0^\circ} \geq 1$ for $i=1,2,\ldots,
l$ and that, for all sequences $(t_i)$,
\begin{equation} \label{upper}
\| \sum_{i=1}^{l} t_i z_i'  \|_{K_0^\circ} \leq \omega \max_i |t_i|,
\end{equation}
where $\omega \leq (4M_1)^{1/2^d} < 4$. Moreover, the length $l$ of the
sequence satisfies
\begin{equation} \label{length}
l \geq \lfloor |\tau|^{1/2^d} \rfloor \geq \lfloor |\tau|^{1/\log_2
  M_1} \rfloor \geq c'' \left(\frac{k} {1 + \log
    \frac{a_1}{a_k}}\right)^{1/2 \log_2 (4M_1)} 
\end{equation}
(note that in our setting we clearly have  $M_1 \geq 1$ and so $d\geq 1$).

\medskip The last step is based on another result from  \cite{Tal}.
\begin{fact} Let $(z_i')_{i=1}^l$ be a sequence in a normed space 
such that $\|z_i'\| \geq 1$ for $i=1,2,\ldots,s$. 
Set $w'= \max \|\sum_{j=1}^s  \ep_iz_i'\|$.
Then there exists a subset $\tau' \subset \{1,2,\ldots , l\}$ 
with $|\tau'| \geq s/8w'$ such that, for any scalars
  $(t_i)$,
$$
\|\sum_{i\in \tau'} t_i z_i \| \geq \frac 12  \max_{i\in \tau'} |t_i | .
$$
\end{fact} 
In our setting,  by (\ref{upper}), 
$m:= |\tau'| > l/32$. On the other hand, by (\ref{length}), the subspace of $Y^*$ spanned by
$z_i', i \in \tau'$, is 8-isomorphic to $\ell_\infty^m$ and hence 
automatically 8-complemented in $Y^*$. The conclusion of
Proposition \ref{ell1} follows then by duality, the only point needing clarification 
being the lower bound on $m$. To elucidate this last issue, we note that the exponent 
$1/\gamma = 1/2\log_2({32 \ell(K_0^\circ)}/{a })$ from the Proposition coincides 
with the lower bound on the exponent 
${1/2 \log_2 (4M_1)} $ in (\ref{length}) given by (\ref{em1}). 
Furthermore, $a_k \geq a$ and 
$$a_1= \|Id : K_0 \to \ell_2^n\| 
= \|Id :  \ell_2^n \to K^\circ_0 \| \leq \ell(K^\circ_0),$$
hence ${a_1}/{a_k} \leq  {\ell(K^\circ_0)}/{a}$ and so, taking again into account the 
form of the lower bound on ${1/2 \log_2 (4M_1)} $ that we are using, 
we conclude that the effect of the 
quantity $1 + \log ({a_1}/{a_k})$ in (\ref{length}) reduces to a multiplicative numerical constant 
(about $0.91$ under the worst case scenario).

\section{Near optimality, and finite dimensional subspaces~of~$L_q$} \label{opt}
The purpose of this section is to substantiate Remark (c), which followed Theorem 
\ref{main} and which asserted that the Theorem as stated can not be essentially 
improved, even if  $X$ varies only over the class of $\ell_q^n$ spaces. To see this, 
denote by  $\gamma_q(Y)$  the factorization constant of $Id_Y$, the identity on $Y$, through 
 an $L_q$-space (i.e., 
 $\gamma_q(Y) := \inf \{\|u\|\,\|v\| \; : \; u:Y\to L_q, \, v : L_q \to Y, \, v\circ u = Id_Y\}$), 
 and similarly $\gamma_q^{(n)}(Y)$ -- the factorization constant of $Id_Y$ through $\ell_q^n$.
 We then have
\begin{fact} \label{factoring} If $q \geq 2$, then

\smallskip 
\noindent  
{\rm 1.} $\gamma_q(\ell_\infty^k) \geq   \,k^{1/q}$

\noindent  
{\rm 2.}   $\gamma_q(\ell_1^k) \geq c \sqrt{k}$ 

\noindent   
{\rm 3.}  $\gamma_q(\ell_2^k) \geq c\min\left\{\sqrt{k }, \sqrt{{q}}\right\}$

\noindent   
{\rm 4.}  $\gamma_q^{(n)}(\ell_2^k) \geq c {\sqrt{k }}/{n^{1/q}}$,

\smallskip\noindent   where $c>0$ is a universal constant.
\end{fact}
\noindent  The ``near optimality" of the statement in the 
Theorem follows now from the fact that if $k$ 
and (large, but not too large) $q >2$ are appropriately related, then all of the quantities in 1.-3. 
must be at least $(\log k)^{1/2}$ (modulo lower order factors; note that Theorem \ref{main} gives an 
upper estimate with exponent 2 in place of 1/2).  
Specifically, if $k$ is sufficiently large and if $q=\log{k}/\log \log {k}$, 
then $k^{1/q} = \log {k}$ and so 
the smallest of the lower bounds is the second expression 
from 3., i.e.,  
$c\sqrt{{q}} = c \sqrt{\log{k}/\log \log {k}}$\,.

\smallskip
The second part of Remark (c) addressed the  ``near optimality" of our estimate 
on the growth of $k_n=c \hskip.2 mm \exp{(\frac 12 \sqrt{\log n}}) $. One possible way 
of stating this assertion more precisely is:
{\em  if for every  $n$-dimensional space $X$ the factorization constant of $Id_{\ell_p^k}$ 
 through $X$ is, for some $p \in \{1, 2, \infty\}$, smaller than $ \exp{( \sqrt{\log k}})$, then 
$k <  \exp{(C    (\log n)^{2/3}})$.}  The argument involves balancing the bounds from 
 1. and from 4. and goes roughly as follows.  
Consider $X=\ell_q^n$, where  $q = \sqrt{\log k}$. 
Then $k^{1/q} = \exp{( \sqrt{\log k}})$, 
which excludes $p=\infty$ (and $p=1$ if $k$ is sufficiently large, which we may assume).
Next, given that  $q = \sqrt{\log k}$, a straightforward calculation shows that 
$k \geq \exp{((4\log n)^{2/3}})$ implies (in fact is equivalent to)
 $n^{1/q} \leq k^{1/4}$ and subsequently implies 
 $c\, \sqrt{k }/{n^{1/q}} \geq c \,k^{1/4} \gg \exp{( \sqrt{\log k}})$.
 This excludes $p=2$.
 If we want to exclude factorization constants smaller than a power of $\log k$ (say,
 $(1+\log{k})^A$, as opposed to $ \exp{( \sqrt{\log k}})$),  
 the argument will be slightly more involved 
 and the resulting restriction on the growth of $(k_n)$ will be -- up to 
 constants depending on $A$ appearing in several places in the exponent -- 
 of the form $\exp{( \sqrt{\log n\, \log \log n}})$.

 \smallskip The argument above may exist in the literature or is a
 folklore; it certainly follows from well-known results and
 methods.  (Indeed, similar considerations might have
motivated various versions of {\em The modified
   problem}; note that it is easy to see that the answer to
 that problem, as stated in the introduction, is affirmative
 if we restrict our attention to spaces $X_n =
 \ell_{q_n}^{m_n}$.)  Similarly, the estimates from Fact
 \ref{factoring} are well known to specialists. In fact,
 the exact values of most (or perhaps even all) quantities
 involved there have been computed.  However, the results
 are spread over the literature and often are not explicitly
 stated. For completeness, we will sketch derivations of
 Fact \ref{factoring} from better known results. (For
 definitions of unexplained concepts and cited facts we
 refer the reader to \cite{PiBook} or \cite{TJ}.)

\smallskip 
\noindent 1.  It is an elementary fact that $d(\ell_\infty^k,\ell_2^k)$, 
the Banach-Mazur distance between $\ell_\infty^k$ and $\ell_2^k$, equals $k^{1/2}$.
A less elementary, but classical estimate (see \cite{Le}) is  that for any 
$k$-dimensional subspace $F \subset L_q$ we have $d(\ell_2^k,F) \leq k^{|1/2-1/q|}$.
Combining these two results we infer that, for any such $F$,  $d(\ell_\infty^k,F) \geq k^{1/q}$.
{\em A fortiori}, $\gamma_q(\ell_\infty^k) \geq k^{1/q}$.

\smallskip 
\noindent 2. By duality, $\gamma_q(\ell_1^k) = \gamma_{q^*}(\ell_\infty^k)$, where
$q^* = q/(q-1) \in [1,2]$ is the dual exponent. Now, the cotype 2 constant of $\ell_\infty^k$ 
is $\sqrt{k}$, while the cotype 2 constants of spaces $L_r$, $1 \leq r \leq 2$ are bounded 
by a universal constant, say  $C$. By the ideal property of the cotype 2 constant 
it follows that, for such $r$,  $ \gamma_{r}(\ell_\infty^k) \geq C^{-1} \sqrt{k}$, and 
the asserted estimate follows.

\smallskip 
\noindent 4. The exact value of $\gamma_\infty(\ell_2^k)$,
the projective constant of $\ell_2^k$,  is well known (\cite{Gr, Ru}), 
in particular we have $\gamma_\infty(\ell_2^k)/\sqrt{k} \in (\sqrt{2/\pi},1]$ for all 
$k \in \mathbb{N}$ (\cite{TJ}, Theorem 32.9(ii); in modern parlance, this is 
a consequence of the ``little" Grothendieck theorem). Consequently, for any $n \in \mathbb{N}$,
$\gamma_\infty^{(n)}(\ell_2^k) \geq \sqrt{2/\pi} \,\sqrt{k}$. This settles the case $q=\infty$, 
and the general case follows since  $d(\ell_\infty^n,\ell_q^n)~=~n^{1/q}$.

\smallskip 
\noindent 3. Again, the estimates for (and even the the exact values of) 
$\gamma_q(\ell_2^k)$ are known to specialists,
but finding them in the literature seems to require combining formulae from 
several sources. First,  
$\gamma_q(\ell_2^k) = n/\pi_q(Id_{\ell_2^n})\pi_{q^*}(Id_{\ell_2^n})$   (\cite{GLR, Re});
this follows from the duality theory for the $\gamma_q$ ideal norm 
(see, e.g., \cite{TJ}, Theorem 13.4) and from symmetries of 
the Hilbert space (cf. \cite{TJ}, \S16). (In fact we need here only the {\em lower} bound on 
$\gamma_q(\ell_2^k)$, which follows just from the duality theory.)
Next, the exact values of, and/or the estimates for $\pi_r(Id_{\ell_2^n})$ 
can be found in \cite{Go, Ga} or in \cite{TJ}, Theorem 10.3. 
And here is a more transparent argument which gives just a sightly weaker estimate
with $\sqrt{q}$ replaced by $\sqrt{q/\log q}$. (This has only minor effect on 
our applications of Fact \ref{factoring}: the lower bound $c \sqrt{\log{k}/\log \log {k}}$ 
becomes $c \sqrt{\log{k}}/\log \log {k}$.)
If $\dim Y = k$, then 
$\gamma_q^{(n)}(Y) \leq 4 \gamma_q(Y)$ for some $n \leq(Ck)^{k}$.
This is because every $k$-dimensional subspace of $L_q$ is contained in a larger 
subspace of dimension $n  \leq (Ck)^{k}$, whose Banach-Mazur distance to $\ell_q^n$ 
is less than (say) 2, and which is 2-complemented in $L_p$ (\cite{PR}).  
Now, if $k \leq q/\log{q}$, then 
$k \log(Ck) \leq q$ (at least for sufficiently large $q$)  and so, for $n$ as above,  
$n^{1/q} \leq \left((Ck)^{k}\right)^{1/q} = \exp(k \log(Ck)/q) \leq e$. 
We now appeal to 4. to deduce that, for all such $k$ and $n$,  
$$
 \gamma_q(\ell_2^k) \geq \frac 14 \gamma_q^{(n)}((\ell_2^k) 
 \geq \frac 14 \sqrt{\frac 2 \pi}\frac {\sqrt{k}} {n^{1/q}} \geq  \frac 1{4e} \sqrt{\frac 2 \pi} {\sqrt{k}} ,
 $$ 
as claimed.  The remaining case $k > q/\log{q}$ follows then from the fact that, for 
fixed $q$, the sequence 
$ \gamma_q(\ell_2^k)$, $k=1,2,\ldots$ is (clearly) nondecreasing.

\smallskip 
{\em Note}: The second argument above would yield the precise version of  3. if we knew that every $k$-dimensional subspace of $L_q$ is contained in a larger subspace whose dimension is (at most) exponential in $k$ and which is, say, 2-isomorphic to $\ell_q^N$ and 2-complemented.
It would be of (independent) interest to clarify this issue, which is relevant to well studied 
``uniform approximation function" of $L_p$-spaces (see \cite{Bo} and its references for 
the background and related results).

\vskip1cm
\small
\address
\end{document}